\documentstyle[12pt]{ioplppt}

\begin{document}
\jl{1}
\title{Simplified boson realization of the $so_q(3)$ subalgebra of 
$u_q(3)$
and matrix elements of $so_q(3)$ quadrupole operators}[Matrix 
elements of
$so_q(3)$ quadrupole operators]

\author{P P Raychev\dag\S, R P Roussev\S, P A Terziev\S,\\
D Bonatsos\ddag\P \ and N Lo Iudice\dag}

\address{\dag\ Dipartimento di Scienze Fisiche, Universit\`a di Napoli
``Federico II''\\
Mostra d' Oltremare, Pad. 19, I-80125 Napoli, Italy}
\address{\S\ Institute for Nuclear Research and Nuclear Energy\\
Bulgarian Academy of Sciences, 72 Tzarigrad Road, BG-1784 Sofia, 
Bulgaria}
\address{\ddag\ European Centre for Theoretical Studies in Nuclear 
Physics
and Related Areas\\
ECT, Villa Tambosi, Strada delle Tabarelle 286, I-38050 Villazzano, 
Trento,
Italy}
\address{\P\ Institute of Nuclear Physics, N.C.S.R. ``Demokritos''\\
GR-15310 Aghia Paraskevi, Attiki, Greece}

\begin{abstract}
A simplified boson realization of the $so_q(3)$ subalgebra of $u_q(3)$
is constructed. A simplified form of the corresponding $so_q(3)$ basis
states is obtained. The reduced matrix elements of a special second-
rank
tensor operator (quadrupole operator) are calculated in the $so_q(3)$ 
basis.
\end{abstract}

\submitted
\maketitle

\section{Introduction}
The construction of chains of subalgebras of a given q-algebra is a
non-trivial problem, since the existence of a chain of subalgebras of 
the
corresponding Lie algebra does not guarantee the existence of the q-
analogue
of this chain. In particular, the $so_q(3)$ subalgebra of $u_q(3)$ has
attracted much attention [1 -- 10], since its classical analogue is a
basic ingredient of several nuclear models, as the Elliott model 
\cite{Ell},
the $su(3)$ limit of the Interacting Boson Model (IBM) \cite{IA} and
the Interacting Vector Boson Model (IVBM) \cite{Ray1}.
The aim of the present paper is to compute the matrix elements of the
$so_q(3)$ quadrupole operator in the $u_q(3) \supset so_q(3)$ basis 
(for the
most symmetric $u_q(3)$-representation). To this purpose we use the 
results
obtained in \cite{J1,Q2}.

In section 2 we introduce a set of modified operators, in terms of 
which the
elements of $so_q(3)$ algebra, i.e. the operators of the q-deformed 
angular
momentum, are expressed in a relatively simple form. In section 3 we 
express
the basis of the q-deformed $so_q(3)\subset u_q(3)$ for the case of 
the most
symmetric representation $[\lambda,0,0]$ of $u_q(3)$. In section 4 we 
also
construct $so_q(3)$ vector operators and in section 5 the reduced 
matrix
elements of a special second-rank tensor operator (quadrupole 
operator) are
calculated in the $so_q(3)$ basis.

\section{Simplified form of the $so_q(3)$ subalgebra of $u_q(3)$}
   In this paper we follow the approach of Refs. \cite{J1,J2}, in
   which a boson realization of the $so_q(3)$ subalgebra of $u_q(3)$ 
in
   terms of q-deformed bosons \cite{Bie,Mac} is constructed.
   The three independent q-deformed boson operators $b_i$ and 
$b_i^\dagger$
   ($i=+,0,-$) satisfy the commutation relations
\begin{equation}\label{eq:s1}
[N_i,b_i^\dagger ] = b_i^\dagger\qquad
[N_i,b_i] = -b_i\qquad
b_i b_i^\dagger -  q^{\pm 1} b_i^\dagger b_i = q^{\mp N_i}
\end{equation}
   where $N_i$ are the corresponding number operators.

   It was shown \cite{J1} that in the Fock space of the totally 
symmetric
   representations $[N,0,0]$ of $u_q(3)$ the angular momentum 
operators i.e.
   the elements of $so_q(3)$ algebra, have the form
\begin{eqnarray}\label{eq:s3}
L_0 &=& N_+ -N_- \nonumber\\
L_+ &=& q^{N_- - \frac{1}{2}N_0}\sqrt{q^{N_+} + q^{-N_+}} b_+^\dagger 
b_0 +
b_0^\dagger b_- q^{N_+ - \frac{1}{2}N_0}\sqrt{q^{N_-} + q^{-N_-
}}\nonumber\\
L_- &=& b_0^\dagger b_ + q^{N_- - \frac{1}{2}N_0}\sqrt{q^{N_+}+q^{-
N_+}} +
q^{N_+ - \frac{1}{2}N_0}\sqrt{q^{N_-} + q^{-N_-}} b_-^\dagger b_0
\end{eqnarray}
   and satisfy the commutation relations
\begin{equation}\label{eq:s4}
[L_0,L_{\pm}] = \pm L_{\pm} \qquad [L_+,L_-] = [2L_0]
\end{equation}
   where the $q$-numbers are defined as $[x] = (q^x-q^{-x})/(q-q^{-
1})$.
   The Casimir operator of $so_q(3)$ can be written in the form 
\cite{Sm}
\begin{eqnarray}\label{eq:s5}
C_2^{(q)} &=& \frac{1}{2}\left\{L_+ L_- + L_- L_+ + [2][L_0]^2
\right\} \nonumber\\
&=& L_- L_+ + [L_0][L_0 + 1] = L_+ L_- + [L_0][L_0 - 1] .
\end{eqnarray}
   In order to rewrite eqs.(\ref{eq:s3}) in a more simplified form, we
   introduce the operators
\begin{eqnarray}\label{eq:s6}
&&B_0 = q^{-\frac{1}{2}N_0}b_0 \qquad\qquad\quad
B_0^\dagger = b_0^\dagger q^{-\frac{1}{2}N_0}\nonumber\\
&&B_i = q^{N_i+\frac{1}{2}} b_i \sqrt{\frac{[2 N_i]}{[N_i]}} \qquad
B_i^\dagger = \sqrt{\frac{[2 N_i]}{[N_i]}} b_i^\dagger 
q^{N_i+\frac{1}{2}}
\qquad i=+,-.
\end{eqnarray}
   These operators satisfy the usual commutation relations
\begin{equation}\label{eq:s7}
[N_i,B^{\dagger}_i] = B_i^{\dagger}\qquad [N_i,B_i] = -B_i .
\end{equation}
   One can check that in the Fock space, spanned on
   the normalized eigenvectors of the excitation number operators
   $N_{+}, N_{0}, N_{-}$, the operators (\ref{eq:s6}) satisfy the 
relations
\begin{eqnarray}\label{eq:s8}
B^{\dagger}_{0} B_{0} = q^{- N_{0}+1}[N_{0}]
&\qquad & B_{0} B^{\dagger}_{0} = q^{-N_{0}}[N_{0}+1]\nonumber\\
B_i^\dagger B_i = q^{2N_i-1}  [2N_i] &\qquad& B_i B_i^\dagger =
q^{2N_i+1} [2N_i+2]\qquad i=+,-
\end{eqnarray}
   from which follow the commutation relations
\begin{equation}\label{eq:s9}
[B_{0},B^{\dagger}_{0}]= q^{-2N_{0}}\qquad
[B_i,B_i^\dagger] = [2] q^{4 N_i+1}\qquad i=+,-.
\end{equation}
   In terms of modified operators (\ref{eq:s6}) the angular momentum
   operators (\ref{eq:s3}) take the simplified form
\begin{eqnarray}\label{eq:s10}
L_{0} &=& N_{+}-N_{-}\nonumber\\
L_{+} &=& q^{-L_{0}+{1\over 2}}B^{\dagger}_{+} B_{0} + q^{L_{0}-
{1\over 2}}
B^{\dagger}_{0} B_{-}\nonumber\\
L_{-} &=& q^{-L_{0}-{1\over 2}}B^{\dagger}_{0} B_{+} + 
q^{L_{0}+{1\over 2}}
B^{\dagger}_{-} B_{0}.
\end{eqnarray}
   It should be noted, however, that these expressions are not 
invariant
   with respect to the replacement $q \leftrightarrow q^{-1}$, which
   restricts us to real $q$.

\section{$so_q(3)$-basis states}
   Using (\ref{eq:s10}) one can check that the normalized highest 
weight
   $so_q(3)$ state $|L L\rangle_q$, which satisfies the conditions
\begin{displaymath}
L_{+}|L L\rangle_q = 0\qquad L_0|L L\rangle_q = L |L L\rangle_q
\qquad\mbox{and}\qquad {}_q{\langle}LL|LL{\rangle}{}_q = 1
\end{displaymath}
   can be written in the form
\begin{equation}\label{eq:b1}
|L L\rangle_q = q^{-{1\over 2}L^2}{(B^{\dagger}_{+})^{L}
\over \sqrt{[2L]!!}}|0\rangle = 
\frac{(b^\dagger_+)^L}{\sqrt{[L]!}}|0\rangle.
\end{equation}
   However these states are {\it not} the most general $so_q(3)$ 
states,
   since they can be multiplied by an arbitrary $so_q(3)$ scalar, 
which
   will not modify the value of $L$. In terms of the modified 
operators one
   can introduce the $so_q(3)$ scalars: \cite{J1,J2}
\begin{eqnarray}\label{eq:b2}
\widetilde S_{+} &=& {1\over [2]} S_{+}  = {1\over [2]}
\left\{(B^{\dagger}_{0})^2 q^{2S_{0}} - {B^{\dagger}_{+} 
B^{\dagger}_{-}}
q^{-2 S_{0} }\right\}\nonumber\\
\widetilde S_{0} &=&  S_{0} = {1\over 2}\left\{N_{+}+N_{0}+N_{-
}+{3\over 2}
\right\} = {1\over 2}\left\{N+{3\over 2} \right\}\nonumber\\
\widetilde S_{-} &=& {1\over [2]} S_{-} = {1\over [2]}
\left\{q^{2 S_{0}}(B_{0})^2 - q^{-2S_{0}}{B_{+} B_{-}}\right\}.
\end{eqnarray}
   These operators satisfy the commutation relations
\begin{equation}\label{eq:b3}
[\widetilde S_{0},\widetilde S_{\pm}] = \pm\widetilde S_{\pm}\qquad
[\widetilde S_{+},\widetilde S_{-}] = -{[2\widetilde S_{0}]}_{q^2}.
\end{equation}
   From (\ref{eq:b3}) it is clear that the set of $so_q(3)$ scalars
   $\widetilde S_{\pm}, \widetilde S_{0}$ close an
   $su_{q^2}(1,1) \sim sp_{q^2}(2,R)$ algebra. Constructing the 
basis, it
   will be simpler to use the scalars $S_{\pm}$, which satisfy the
   commutation relations
\begin{equation}\label{eq:b4}
[S_{-},S_{+}] = [2]^2 {[2S_{0}]}_{q^2} = [2][2N+3]
\qquad (S_{+})^{\dagger} = S_{-}.
\end{equation}

   Therefore, the so$_q$(3) states, characterized by an angular 
momentum $L$
   and its projection $M = L$, which belong to the most symmetric
   $[\lambda, 0, 0]$ irreducible representation of $u_q(3)$ can be 
written in
   the form
\begin{equation}\label{eq:b5}
\left\vert\begin{array}{cc}{\lambda}&{}\\{L}&{L}\end{array}\right\rangle
= {1\over N_{\lambda L}}(S_{+})^{{1\over 2}(\lambda-L)}|L L\rangle_q
\end{equation}
   where $L = \lambda,\lambda-2,\ldots,0$ or $1$ and $|L L\rangle_q$
   is a notation for the states (\ref{eq:b1}). The normalization 
constant
   $N_{\lambda L}$ is determined from the condition
\begin{displaymath}
{{}\atop{}}_{q}{\left\langle{\begin{array}{cc}{\lambda}&{}\\{L}&{L}
\end{array}}\right.}{{\left\vert\begin{array}{cc}{\lambda}&{}\\{L}&{L}
\end{array}\right\rangle}\!{{}\atop {}}_{q}} = 1.
\end{displaymath}
   Using the relations
\begin{eqnarray*}
[S_{-},(S_{+})^k] = [2k] (S_{+})^{k-1} [2N+2k+1]\\
S_{-}|L L\rangle_q = 0\\
(S_{-})^k (S_{+})^k |L L\rangle_q = {[2k]!![2L+2k+1]!! \over [2L+1]!!}
|L L\rangle_q \qquad\quad k={1\over 2}(\lambda-L)
\end{eqnarray*}
   the final result is
\begin{equation}\label{eq:b6}
N_{\lambda L} = \sqrt{[\lambda-L]!![\lambda+L+1]!! \over [2L+1]!!}.
\end{equation}

Now the states with an arbitrary projection of the momentum $M$ can be
obtained by successive application of $L_{-}$ on the states 
(\ref{eq:b5}),
i.e.
\begin{eqnarray}\label{eq:b7}
\fl  {\left\vert\begin{array}{cc}{\lambda}&{}\\{L}&{M}\end{array}
\right\rangle}_{q} = \sqrt{[L+M]! \over [2L]![L-M]!}
(L_{-})^{L-M} 
{\left\vert\begin{array}{cc}{\lambda}&{}\\{L}&{L}\end{array}
\right\rangle}_{q} \nonumber\\
\lo{=} \frac{q^{-{1\over 2}L^2}}{N_{\lambda L}}\sqrt{[L+M]! \over 
[2L]![L-M]!}
(S_{+})^{{1\over 2}(\lambda-L)} (L_{-})^{L-M}
{(B^{\dagger}_{+})^L \over \sqrt{[2L]!!}}|0\rangle.
\end{eqnarray}
   In order to find an explicit expression for the states 
(\ref{eq:b7}) in
   form of a polynomial in terms of the operators $B^\dagger_i$
   we shall make use of the auxiliary formula
\begin{eqnarray}\label{eq:b9}
\fl  L_{-} {(B^{\dagger}_{+})^x \over [2x]!!}{(B^{\dagger}_{0})^y 
\over [y]!}
{(B^{\dagger}_{-})^z \over [2z]!!} |0\rangle
=& q^{x-y-z+{1\over 2}}[2z+2]{(B^{\dagger}_{+})^x \over [2x]!!}
{(B^{\dagger}_{0})^{y-1}
\over [y-1]!}{(B^{\dagger}_{-})^{z+1} \over [2z+2]!!} |0\rangle 
\nonumber\\
& + q^{x+z-{1\over 2}}[y+1]{( B^{\dagger}_{+} )^{x-1} \over [2x-2]!!}
{(B^{\dagger}_{0})^{y+1} \over [y+1]!}{(B^{\dagger}_{-})^z \over 
[2z]!!}
|0\rangle.
\end{eqnarray}
   where $x \geq 1, y \geq 1$ and $z \geq 0$. Using (\ref{eq:b9}) one 
can
   prove by induction in $m \geq 0$ that the following relation holds
\begin{equation}\label{eq:b9i}
\fl  (L_{-})^m {(B^{\dagger}_{+})^L \over [2L]!!}|0\rangle =
q^{{1 \over 2} m(2L-m)}[m]!
\sum_{p}{(B^{\dagger}_{+})^p \over [2p]!!}{(B^{\dagger}_{0})^{2L-m-2p}
\over [2L-m-2p]!}{(B^{\dagger}_{-})^{m-L+p} \over [2m-
2L+2p]!!}|0\rangle
\end{equation}
   where the summation index $p$ runs over these values, for which
   all exponents of the operators $B_{i}^{\dagger}$
   are not negative. Replacing $m = L - M$ in (\ref{eq:b9i}) we obtain
\begin{equation}\label{eq:b10}
\fl  {(L_{-})^{L-M} \over [L-M]!}{(B^{\dagger}_{+})^L \over 
[2L]!!}|0\rangle =
q^{{1\over 2}(L^2-M^2)}
{\sum_{p=max(0,M)}^{\left\lfloor (L+M)/2 \right\rfloor}}
{(B^{\dagger}_{+})^p \over [2p]!!}{(B^{\dagger}_{0})^{L+M-2p}
\over [L+M-2p]!}{(B^{\dagger}_{-})^{p-M} \over [2p-2M]!!}|0\rangle.
\end{equation}
   After combining (\ref{eq:b6}), (\ref{eq:b7}) and (\ref{eq:b10}) we 
obtain
   the following expression for the basis states
\begin{eqnarray}\label{eq:b11}
\fl  {\left\vert\begin{array}{cc}{\lambda}&{}\\{L}&{M}\end{array}
\right\rangle}_{q}
= q^{-{1\over 2}M^2}\sqrt{[L+M]![L-M]![2L+1] \over
[\lambda-L]!![\lambda+L+1]!!} \nonumber\\
 \times (S_{+})^{{1\over 2}(\lambda-L)}
{\sum_{p=max(0,M)}^{\left\lfloor (L+M)/2 \right\rfloor}}
{(B^{\dagger}_{+})^p \over
[2p]!!}{(B^{\dagger}_{0})^{L+M-2p} \over [L+M-2p]!}
{(B^{\dagger}_{-})^{p-M} \over [2p-2M]!!}|0\rangle.
\end{eqnarray}

   In order to rewrite the basis states (\ref{eq:b11}) in a polynomial
   form, by expanding the power of $S_{+}$, one can use the q-binomial
   theorem \cite{Nom}, according to which, if the elements $X$ and $Y$
   satisfy the condition $Y X = q X Y$ then
\begin{equation}\label{eq:b13}
{(X - Y)}^k = \sum_{t=0}^k (-1)^t q^{{1\over 2}t(k-t)}
{\left[ \begin{array}{c}{k}\\{t}\end{array} \right] }_{q^{1\over 2}}
X^{k-t} Y^t.
\end{equation}
   In the present case we have
\begin{displaymath}
S_{+} = \underbrace{(B^{\dagger}_{0})^2 q^{2S_0}}_{X} -
\underbrace{B^{\dagger}_{+} B^{\dagger}_{-} q^{-2S_0}}_{Y}\qquad\quad
Y  X = q^{-4} X Y.
\end{displaymath}
   Therefore for the power of $S_{+}$ we obtain
\begin{equation}\label{eq:b14}
\fl  (S_{+})^k = \sum_{t=0}^k (-1)^t q^{-2t(k-t)}
{\left[\begin{array}{c}{k}\\{t}\end{array} \right] }_{q^2}
\left\{ (B^{\dagger}_{0})^2 q^{2S_{0}} \right\}^{k-t}
\left\{ B^{\dagger}_{+} B^{\dagger}_{-} q^{-2S_{0}} \right\}^t
\end{equation}
   where
\begin{displaymath}
\left[ \begin{array}{c} k \\ t \end{array} \right]_{q^2} =
{[k]_{q^2}! \over [t]_{q^2}! [k-t]_{q^2}!} =
{[2k]!! \over [2t]!![2k-2t]!!}\qquad 2S_{0} = N + {3\over 2}
\end{displaymath}
   and grouping the terms with $q^N$ we have
\begin{equation}\label{eq:b15}
\fl  (S_{+})^k = q^{k(k+{1\over 2})} [2k]!! \sum_{t=0}^k
\frac{(-1)^t q^{-(2k+1)t}}{[2t]!![2k-2t]!!} (B^{\dagger}_{+})^t
(B^{\dagger}_{0})^{2(k-t)} (B^{\dagger}_{-})^t q^{(k-2t)N}.
\end{equation}
   Combining eqs.(\ref{eq:b11}) and (\ref{eq:b15}) for
   $k = {1\over 2}(\lambda-L)$ the basis states (\ref{eq:b7})
   can be written in the form \cite{J1,J2,J4}
\begin{eqnarray}\label{eq:b16}
\fl  {\left\vert\begin{array}{cc}{\lambda}&{}\\{L}&{M}\end{array}
\right\rangle}_{q} = q^{{1\over 4}(\lambda-L)(\lambda+L+1) -
{1\over 2}{M^2}} \sqrt{[L+M]![L-M]![\lambda-L]!![2L+1] \over
[\lambda+L+1]!!}  \nonumber\\
\fl\quad  \times\sum_{t=0}^{(\lambda-L)/2}
{\sum_{p=max(0,M)}^{\left\lfloor (L+M)/2 \right\rfloor}}
{{(-1)^t}q^{-(\lambda+L+1)t} \over [2t]!![\lambda-L-2t]!!}
{(B^{\dagger}_{+})^{p+t} \over [2p]!!}{(B^{\dagger}_{0})^{\lambda+M-
2p-2t}
\over [L+M-2p]!}{(B^{\dagger}_{-})^{p+t-M} \over [2p-2M]!!}|0\rangle.
\end{eqnarray}

\section{Vector operators}
   The $so_q(3)$ tensor operators must satisfy the commutation 
relations,
   which directly follow from the expression for the adjoint action 
of the
   corresponding algebra \cite{Bie2,Sm,Nom}. By definition, the
   irreducible tensor operator $T^j_m$ of rank $j$ according to 
$so_q(3)$
   satisfies the commutation relations
\begin{eqnarray}\label{eq:t8}
[L_0,T^j_m] = m\ T^j_m\nonumber\\{}
[L_\pm,T^j_m]_{q^m} q^{L_0} =
\sqrt{[j \mp m][j \pm m + 1]}\ T^j_{m \pm 1}.
\end{eqnarray}
   The generalization of the Wigner-Eckart theorem to the case of the
   algebra $so_q(3)$ is
\begin{equation}\label{eq:t8a}
\langle\alpha',L' M'|T^j_m|\alpha,L M\rangle = (-1)^{2j}\;
\frac{_qC_{LM,jm}^{L'M'}}{\sqrt{[2L'+1]}}
\langle\alpha',L'\|T^j\|\alpha,L\rangle
\end{equation}
   where $|\alpha,L M\rangle$ are orthonormalized basis vectors of the
   irreducible representation $_qD^L$ of the algebra $so_q(3)$ and
   $_qC_{L_1 M_1,L_2 M_2}^{LM}$ are the Clebsch-Gordan coefficients
   \cite{Sm,Nom} of the same algebra. It should be noted that, the 
operator
\begin{equation}\label{eq:t9}
R^j_m = (-1)^m q^{-m} (T^j_{-m})^\dagger
\end{equation}
   where the superscript $^{\dagger}$ denotes Hermitian conjugation,
   transforms in the same way (\ref{eq:t8}) as the tensor operator
   $T^j_m$, i.e. it also is an irreducible $so_q(3)$ tensor operator
   of rank $j$.

   In order to construct irreducible $so_q(3)$ vector operators
   $T^\dagger_m$ and $\widetilde{T}_m$ we start from the observation
\begin{equation}\label{eq:v1}
[L_0,B^\dagger_+] = B^\dagger_+\nonumber
\end{equation}
   and suppose that the highest weight component of the vector 
operator
   $T^\dagger_m$ is
\begin{equation}\label{eq:v2}
T^{\dagger}_{+1} = \omega B^{\dagger}_{+} q^{\alpha N_{+} + \beta 
N_{0} +
\gamma N_{-} + \delta}
\end{equation}
   where $\alpha, \beta, \gamma, \delta$ and $\omega$ are real 
constants to be
   determined. As irreducible first-rank $so_q(3)$ tensor operator,
   $T^{\dagger}_{m}\ (m=0,\pm1)$ must satisfy the relations
\begin{eqnarray}\label{eq:v3}
[L_0,T^{\dagger}_{m}] = m\,T^{\dagger}_{m}\nonumber\\{}
[L_{\pm},T^{\dagger}_{m}]_{q^m} q^{L_0} =
\sqrt{[1 \mp m][2 \pm m]} T^{\dagger}_{m\pm1}.
\end{eqnarray}
   The same relations hold for the operators
\begin{equation}\label{eq:v4}
\widetilde{T}_{m} = (-1)^m q^{-m} (T^{\dagger}_{-m})^{\dagger} =
(-1)^m q^{-m} T_{-m}
\end{equation}
   where $(T^{\dagger}_{m})^{\dagger} = T_{m}$ and $^{\dagger}$
   denotes Hermitian conjugation. According to eq.(\ref{eq:v3}), the 
condition
\begin{equation}\label{eq:v5}
{[L_{+},T^{\dagger}_{+1}]}_{q} = 0
\end{equation}
   is satisfied, if $\alpha + 2 = \beta = \gamma$ for any real 
constants
   $\omega$ and $\delta$, and the operator $T^{\dagger}_{+1}$ can be 
written as
\begin{equation}\label{eq:v6}
T^{\dagger}_{+1} = \omega B^{\dagger}_{+} q^{-2 N_{+} + \beta N + 
\delta}.
\end{equation}
   Further by the action of $L_-$ we get all other components of 
$T^\dagger_m$
\begin{eqnarray}\label{eq:v8}
T^{\dagger}_{0} &= \omega \sqrt{[2]} B^{\dagger}_{0} q^{-2 N_{+} + 
\beta N +
\delta + {1\over 2}}\nonumber\\
T^{\dagger}_{-1} &= \omega \left\{B^{\dagger}_{-} q^{2N_{+} + (\beta-
2)N +
\delta} - (q - q^{-1}) B_{+} (B^{\dagger}_{0})^2 q^{-2 N_{+} +
\beta N + \delta +2} \right\}.
\end{eqnarray}
   One can check that the condition
\begin{equation}\label{eq:v11}
{[L_{-},T^{\dagger}_{-1}]}_{q^{-1}} = 0
\end{equation}
   holds for any values of the parameters $\beta, \delta$ and 
$\omega$. From
   these expressions it is clear that $T^\dagger_m\ (m=0,\pm1)$ is a 
vector
   operator according to $so_q(3)$. The components of the 
corresponding
   conjugated vector operator $\widetilde{T}_{m}\ (m=0,\pm1)$ given by
   (\ref{eq:v4}) are
\begin{eqnarray}\label{eq:v12}
\widetilde{T}_{+1} &= -\omega \left\{q^{2N_{+} + (\beta-2)N + \delta -
 1}
B_{-} - (q - q^{-1}) q^{-2 N_{+} + \beta N + \delta + 1} 
B^{\dagger}_{+}
(B_{0})^2 \right\}\nonumber\\
\widetilde{T}_{0} &= \omega \sqrt{[2]} q^{-2N_{+} + \beta N + \delta +
{1\over 2}} B_{0}\nonumber\\
\widetilde{T}_{-1} &= -\omega q^{-2N_{+} + \beta N + \delta + 1} 
B_{+}.
\end{eqnarray}
   Using the vector operators $T^{\dagger}_{m}$ and 
$\widetilde{T}_{m}$
   one can construct the coupled operators \cite{Q2,Sm}
\begin{equation}\label{eq:v13}
A_{M}^{L} = {[T^{\dagger} \otimes \widetilde{T}]}^{L}_{M} =
\sum_{m,n}{}_{q^{-1}}C_{1m,1n}^{LM} T^{\dagger}_{m} \widetilde{T}_{n}
\qquad L=0,1,2.
\end{equation}
   Actually we use a particular case of a product of two irreducible 
tensor
   operators acting on a single vector \cite{Sm}. If $T^{\dagger}_{m}$
   and $\widetilde{T}_{m}$ are vector operators according to 
$so_q(3)$ then
   the operators (\ref{eq:v13}) are irreducible tensors of rank 
$L=0,1,2$
   according to the same algebra. Their Hermitian conjugates are
\begin{equation}\label{eq:v14}
(A^{L}_{M})^{\dagger} = (-1)^M q^{-M} A^{L}_{-M}.
\end{equation}
   In order to determine the parameters $\beta, \delta$ and $\omega$ 
we shall
   take into account that from the generators $L_{+}, L_{0}, L_{-}$ 
of the
   algebra $so_q(3)$ one can construct a first-rank tensor $J^{1}$
   \ \cite{Sm,Nom} according to this algebra as
\begin{eqnarray}
J^{1}_{\pm1} &=& \mp{1\over \sqrt{[2]}} q^{- L_{0} } 
L_{\pm}\label{eq:v15}\\
J^{1}_{0} &=& {1\over [2]}\left\{q L_{+} L_{-} - q^{-1} L_{-} L_{+} 
\right\}
= {1\over [2]}\left\{ q [2L_{0}] + (q - q^{-1}) L_{-} L_{+} \right\}
\label{eq:v15b}\\
&=& {1\over [2]}\left\{ q [2L_{0}] + (q - q^{-1})
\left( C_{2}^{(q)} - [L_{0}][L_{0} + 1] \right) \right\}\nonumber
\end{eqnarray}
   where $C_{2}^{(q)}$ is the second order Casimir operator 
(\ref{eq:s5}) of
   $so_q(3)$. After imposing the condition
\begin{equation}\label{eq:v16}
J_{M}^{1} = -\sqrt{[4]\over [2]} A^{1}_{M}\qquad M=0,\pm1
\end{equation}
   where $A^{1}_{M}$ is a first-rank tensor (\ref{eq:v13}) and 
$J_{M}^{1}$ is
   also a first-rank tensor (\ref{eq:v15}) we obtain
\begin{equation}\label{eq:v18}
\omega = {1\over \sqrt{[2]}}\qquad \beta = 1\qquad \delta = -{1\over 
2}
\qquad \alpha + 2 = \beta = \gamma.
\end{equation}
   The final expressions for the components of the vector operator
   $T^{\dagger}_{m}$ are
\begin{eqnarray}\label{eq:v19a}
&&T^{\dagger}_{+1} = {1\over \sqrt{[2]}} B^{\dagger}_{+}
q^{-2N_{+} + N - {1\over 2}}\nonumber\\
&&T^{\dagger}_{0} = B^{\dagger}_{0}  q^{-2N_{+} + N}\nonumber\\
&&T^{\dagger}_{-1} = {1\over \sqrt{[2]}}\left\{ B^{\dagger}_{-}
q^{2 N_{+} - N - {1\over 2}} - (q - q^{-1}) B_{+} 
{(B^{\dagger}_{0})}^2
q^{-2 N_{+} + N + {3\over 2}}\right\}.
\end{eqnarray}
   The corresponding expressions for $\widetilde{T}_m$ can be 
obtained from
   (\ref{eq:v4})
\begin{eqnarray}\label{eq:v19b}
&&\widetilde{T}_{+1} = -{1\over \sqrt{[2]}}\left\{
q^{2 N_{+} - N - {3\over 2}} B_{-} - (q-q^{-1}) q^{-2N_{+} + N+ 
{1\over 2}}
B^{\dagger}_{+} {(B_{0})}^2 \right\}\nonumber\\
&&\widetilde{T}_{0} = q^{-2N_{+} + N} B_{0}\nonumber\\
&&\widetilde{T}_{-1} = -{1\over \sqrt{[2]}}q^{-2N_{+}+N+{1\over 2}} 
B_{+} .
\end{eqnarray}
   In terms of $T^\dagger_m$ and $\widetilde{T}_m$ one can construct 
the
   scalars
\begin{eqnarray}\label{eq:v20}
R_{+}&=&-\sqrt{[3]} [T^{\dagger} \otimes T^{\dagger}]^{0}_{0} =
-\sqrt{[3]} \sum_{m,n}{}_{q^{-1}}C_{1m,1n}^{00} T^{\dagger}_{m}
T^{\dagger}_{n}\nonumber\\
R_{-}&=&-\sqrt{[3]} [\widetilde{T} \otimes \widetilde{T} ]^{0}_{0} =
-\sqrt{[3]} \sum_{m,n}{}_{q^{-1}}C_{1m,1n}^{00} \widetilde{T}_{m}
\widetilde{T}_{n}
\end{eqnarray}
   and can easily check that
\begin{equation}\label{eq:v21}
R_{+} = S_+ q^{2S_0} \qquad R_- = q^{2S_0} S_- \qquad\qquad 
(R_+)^\dagger = R_-
\end{equation}
   where $S_+, S_0, S_-$ are given by (\ref{eq:b2}).

   Now using the relations
\begin{eqnarray*}
[T^\dagger_{+1},T^\dagger_0]_{q^2} = 0 && \\{}
[T_{+1},T^\dagger_{+1}]_{q^{-2}} = q^{2N} &\quad\mbox{or}\quad&
[\widetilde{T}_{-1},T^\dagger_{+1}]_{q^{-2}} = -q^{2N+1}
\end{eqnarray*}
 which follow from the explicit form (\ref{eq:v19a},\ref{eq:v19b}) of 
these
 operators, one can compute by successive application of $L_-$ the 
commutation
 relations between $\widetilde{T}_m$ and $T^\dagger_m$
\begin{eqnarray}
[T^\dagger_{+1},T^\dagger_0]_{q^2} = 0 &\qquad&
[\widetilde{T}_0,\widetilde{T}_{-1}]_{q^2} = 0 \nonumber\\{}
[T^\dagger_0,T^\dagger_{-1}]_{q^2} = 0 &\qquad&
[\widetilde{T}_{+1},\widetilde{T}_0]_{q^2} = 0 \label{eq:v22a}\\{}
[T^\dagger_{+1},T^\dagger_{-1}] = (q-q^{-1}) (T^\dagger_0)^2 &\qquad&
[\widetilde{T}_{+1},\widetilde{T}_{-1}] = (q-q^{-1}) 
(\widetilde{T}_0)^2
\nonumber
\end{eqnarray}
   and
\numparts
\begin{eqnarray}
[\widetilde{T}_0,T^\dagger_{+1}] = 0 &\qquad&
[\widetilde{T}_{-1},T^\dagger_0] = 0 \nonumber\\{}
[\widetilde{T}_{+1},T^\dagger_{+1}]_{q^2} = 0 &\qquad&
[\widetilde{T}_{-1},T^\dagger_{-1}]_{q^2} = 0 \label{eq:v22b}\\{}
[\widetilde{T}_{+1},T^\dagger_0] = (q^2-q^{-2})
T^\dagger_{+1}\widetilde{T}_0 &\qquad&
[\widetilde{T}_0,T^\dagger_{-1}] = (q^2-q^{-2})
T^\dagger_0\widetilde{T}_{-1} \nonumber
\end{eqnarray}
\begin{eqnarray}
[\widetilde{T}_{-1},T^\dagger_{+1}]_{q^{-2}} = -q^{2N+1} \nonumber\\{}
[\widetilde{T}_0,T^\dagger_0] = q^{2N} + q^{-1}(q^2-q^{-2})
T^\dagger_{+1}\widetilde{T}_{-1} \label{eq:v22c}\\{}
[\widetilde{T}_{+1},T^\dagger_{-1}]_{q^{-2}} = -q^{2N-1} + q^{-1}(q^2-
q^{-2})
\left\{T^\dagger_0\widetilde{T}_0 + (q-q^{-1}) 
T^\dagger_{+1}\widetilde{T}_{-1}
\right\} \nonumber
\end{eqnarray}
\endnumparts
which are similar to the results obtained by Quesne in \cite{Q2}.

\section{Matrix elements of the quadrupole operator}
   The aim of this section is the calculation of the reduced matrix 
elements
   of the q-deformed quadrupole operator $Q^2$, namely, the quantities
\begin{displaymath}
\langle\lambda,L+2\|Q^2\|\lambda,L\rangle\qquad\mbox{and}\qquad
\langle\lambda,L\|Q^2\|\lambda,L\rangle
\end{displaymath}
   where the operator $Q^2$ is
\begin{equation}\label{eq:q0}
Q^2_M = \sqrt{\frac{[3][4]}{[2]}} A^2_M, \quad
A^2_M = {[T^{\dagger} \otimes \widetilde{T}]}^{2}_{M} =
\sum_{m,n}{}_{q^{-1}}C_{1m,1n}^{2M} T^{\dagger}_{m} \widetilde{T}_{n}.
\end{equation}
   In eq.(\ref{eq:q0}) the factor has been chosen to agree with the 
usual
   convention in the classical case when $q\rightarrow 1$. Other 
reduced
   matrix elements do not occur, since in the most symmetric 
representation of
   $u_q(3)$ only states with equal parity of $\lambda$ and $L$ exist. 
From
   the tensor structure of the operator $A^2$ for its zero component 
we have
\begin{equation}\label{eq:q1}
A^2_0 {\left|\begin{array}{cc}{\lambda}&{}\\{L}&{L}\end{array}
\right\rangle}_q =
{\bf a}\,{\left|\begin{array}{cc} {\lambda}&{}\\{L+2}&{L}
\end{array}\right\rangle}_q +
{\bf b}\,{\left|\begin{array}{cc}{\lambda}&{}\\{L}&{L}\end{array}
\right\rangle}_q
\end{equation}
   and it is clear that the coefficients ${\bf a, b}$ determine the
   reduced matrix elements of the tensor $A^2$.

   The highest weight vector of the basis states (\ref{eq:b5}) can be 
expressed
   in terms of vector operators (\ref{eq:v19a}) as follows
\begin{eqnarray}\label{eq:q6}
{\left|\begin{array}{cc}{\lambda}&{}\\{L}&{L}\end{array}\right\rangle}
\frac{(S_+)^k}{N_{\lambda 
L}}\frac{(T^\dagger_{+1})^L}{\sqrt{[L]_{q^2}!}}
|0\rangle
\end{eqnarray}
   where $k = \frac{1}{2}(\lambda - L)$  and the normalization 
constant
   $N_{\lambda L}$ is determined in (\ref{eq:b6}). Therefore
\begin{eqnarray}\label{eq:q9}
A^2_0 
{\left|\begin{array}{cc}\lambda&\\L&L\end{array}\right\rangle}_q &=&
A^2_0 \frac{(S_+)^k}{N_{\lambda 
L}}\frac{(T^\dagger_{+1})^L}{\sqrt{[L]_{q^2}!}}
|0\rangle\nonumber\\
&=& \frac{1}{N_{\lambda L}\sqrt{[L]_{q^2}!}}\left\{ (S_+)^k A^2_0 +
[A^2_0,(S_+)^k] \right\} (T^\dagger_{+1})^L |0\rangle.
\end{eqnarray}
   Now, in order to calculate the action of $A^2_0$ on the highest 
weight
   vector (\ref{eq:q6}), we use the identities
\begin{eqnarray}\label{eq:q8}
&&[A^2_2,(S_+)^k] = q^{2k-2} [2k] (S_+)^{k-1} (T^\dagger_{+1})^2 
q^{2S_0}
\nonumber\\{}\label{eq:q7}
&&[A^2_1,(S_+)^k] = \sqrt{\frac{[4]}{[2]}} q^{2k-1} [2k] (S_+)^{k-1}
T^\dagger_0 T^\dagger_{+1} q^{2S_0}\nonumber\\{}
&&[A^2_0,(S_+)^k] = \sqrt{\frac{[4]}{[3][2]}} q^{2k} [2k] (S_+)^{k-1}
\left\{ S_+ q^{2S_0+1} + [3] T^\dagger_{-1} T^\dagger_{+1} \right\} 
q^{2S_0}
\end{eqnarray}
   which can be obtained from 
eqs.(\ref{eq:v22a},\ref{eq:v22b},\ref{eq:v22c})
   and successive application of the operator $L_-$.
   Finally using the relations
\begin{eqnarray}
A^2_0 (T^\dagger_{+1})^L |0\rangle &=& -\sqrt{\frac{[2]}{[3][4]}}
\frac{q^3 [2L]}{[2]} (T^\dagger_{+1})^L |0\rangle \label{eq:q10}\\
T^\dagger_{-1}(T^\dagger_{+1})^{L+1}|0\rangle &=& \frac{q^{-2L-
2}}{[2L+4][2L+3]}
(L_-)^2(T^\dagger_{+1})^{L+2}|0\rangle \nonumber\\
&&-\,\frac{q^{L+\frac{5}{2}}[2L+2]}{[2][2L+3]} S_+ (T^\dagger_{+1})^L 
|0\rangle
\label{eq:q11}
\end{eqnarray}
   we obtain the expressions for the coefficients ${\bf a, b}$ in
   the expansion (\ref{eq:q1})
\begin{equation}\label{eq:q12a}
{\bf a} = \frac{q^{\lambda-2L-\frac{1}{2}}}{[2L+3]} 
\sqrt{\frac{[3][4]}{[2]}}
\sqrt{\frac{[\lambda-L][\lambda+L+3][2L+2]}{[2][2L+5]}}
\end{equation}
   and
\begin{equation}\label{eq:q12b}
{\bf b} = -\frac{q^{\lambda+{5\over 2}}[2L]}{[2][2L+3]}
\sqrt{\frac{[2]}{[3][4]}} \left\{ q^{L-\frac{1}{2}}[\lambda-L] +
q^{-L+\frac{1}{2}}[\lambda+L+3] \right\}.
\end{equation}

   From eqs.(\ref{eq:q12a},\ref{eq:q12b}) and the Wigner-Eckart 
theorem
   (\ref{eq:t8a}) one immediately obtains the reduced matrix elements 
of the
   quadrupole operator $Q^2$ defined by (\ref{eq:q0})
\begin{eqnarray}
\fl  \langle{\lambda,L+2}\|Q^{2}\|{\lambda,L}\rangle
= {q^{\lambda-{1\over 2}} \over [2]} \sqrt{[3][4] \over [2]}
\sqrt{[\lambda-L][\lambda+L+3][2L+4][2L+2] \over 
[2L+3]}\label{eq:q14a}\\
\fl  \langle{\lambda,L}\|Q^{2}\|{\lambda,L}\rangle
= -{q^{\lambda-{1\over 2}}\over [2]} \sqrt{[2L][2L+1][2L+2] \over
[2L-1][2L+3]}\nonumber\\
\times\left\{ q^{L-{1\over 2}}[\lambda-L] +
q^{-L+{1\over 2}}[\lambda+L+3] \right\} \label{eq:q14b}
\end{eqnarray}
   which is in agreement with the classical case when $q\rightarrow 
1$.
   Taking into account the Wigner-Eckart theorem and the symmetry 
properties
   of the q-deformed Clebsch-Gordan coefficients it can be shown also 
that
   the reduced matrix element (\ref{eq:q14a}) of the q-deformed 
quadrupole
   operator (\ref{eq:q0}) has the following symmetry property
\begin{equation}\label{eq:q15}
\langle\lambda,L+2\|Q^2\|\lambda,L\rangle =
\langle\lambda,L\|Q^2\|\lambda,L+2\rangle .
\end{equation}
   For small values of the deformation parameter $\tau$
   \ ($q={\rm e}^\tau$, $\tau$-real) the reduced matrix elements
   (\ref{eq:q14a},\ref{eq:q14b}) can be represented in  Taylor 
expansions
\begin{eqnarray}\label{eq:q16a}
\fl  \langle\lambda,L+2\|Q^2\|\lambda,L\rangle =
\sqrt{\frac{6(\lambda-L)(\lambda+L+3)(L+2)(L+1)}{2L+3}}\nonumber\\
\fl\qquad  \times\left\{1 + \left(\lambda-\frac{1}{2}\right)\tau +
\left(\frac{2}{3}\lambda^2 + 
\frac{1}{2}L^2+\frac{3}{2}L+\frac{65}{24}\right)
\tau^2 + {\rm O}(\tau^3)\right\}
\end{eqnarray}
\begin{eqnarray}\label{eq:q16b}
\fl  \langle\lambda,L\|Q^2\|\lambda,L\rangle =
-(2\lambda+3)\sqrt{\frac{L(L+1)(2L+1)}{(2L-1)(2L+3)}}
\biggl\{1+\frac{2}{2\lambda+3}\{\lambda(\lambda+1)-L(L+1)\}\tau
\nonumber\\
\fl\qquad  + \frac{1}{3(2\lambda+3)}\{(2\lambda+15)L(L+1)+(2\lambda+1)
(2\lambda^2+2\lambda+3)\}\tau^2 + {\rm O}(\tau^3)\biggr\}
\end{eqnarray}
   and if $\tau=0$ one obtains the classical expressions for the 
corresponding
   reduced matrix elements.

\section{Conclusion}
In the present paper we gave an approach for the construction of
irreducible tensor operators in the case of the q-deformed chain
$u_q(3)\supset so_q(3)$ for the most symmetric representations
$[\lambda,0,0]$ of the $u_q(3)$ algebra. In this way we have 
calculated the
reduced matrix elements (\ref{eq:q14a},\ref{eq:q14b}) of the q-
deformed
quadrupole operator (\ref{eq:q0}).

>From physical point of view of great interest are the $E2$-transition
probabilities ($B[E2]$-factors), which in the classical case are 
expressed
by means of the reduced matrix elements of the $u(3)$-quadrupole 
operator.
It can be shown that the $B[E2]$-factors corresponding to the chain
$u_q(3)\supset so_q(3)$ are of the form
\begin{equation}\label{eq:q37}
B[E2; (\lambda,L+2) \rightarrow (\lambda,L)]_q =
\frac{1}{[2L+5]}{|\langle\lambda,L+2\|Q^2\|\lambda,L\rangle|}^2
\end{equation}
where $Q^2$ is the q-deformed quadrupole operator (\ref{eq:q0}). 
Likewise
the reduced matrix element (\ref{eq:q14b}) is related to the 
deformation of
the physical system in the state with angular momentum $L$.

It should be noted that the results obtained here, strictly speaking, 
are
valid only for real values of the deformation parameter $q$. On the 
other
hand the comparison of the experimental data with the predictions of 
a number
of physical models, \cite{Ray2,Kib}
based on the q-deformed $su_q(2)$ algebra, shows that one
can achieve a good agreement between theory and experiment only if 
$q$ is a
pure phase ($q = {\rm e}^{i\tau}$). Nevertheless, we suppose, however,
that these quadrupole operators describe some q-deformed excitations
(q-deformed phonons) and the obtained results are a necessary step to
further investigations.

\ack
One of the authors (DB) has been supported by the EU under contract
ERBCHBGCT 930467. Another authors (PPR) and (RPR) have been supported 
by
the Bulgarian Ministry of Science and Education under contracts 
$\Phi$-415
and $\Phi$-547.

\end{document}